\documentclass[12pt]{amsart}
\usepackage[margin=3.5cm]{geometry}
\usepackage{amssymb}
\usepackage{relsize}
\usepackage{url}
\numberwithin{equation}{section}

\newcommand{\di}{d_i}

\newcommand\numeq[1]%
  {\stackrel{\scriptscriptstyle(\mkern-1.5mu#1\mkern-1.5mu)}{=}}

\newcounter{relctr} 
\everydisplay\expandafter{\the\everydisplay\setcounter{relctr}{0}} 

\AtBeginDocument{} 

\newcommand{\jo}{{j_1}}
\newcommand{\jtw}{{j_2}}

\usepackage{tikz-cd}
\usepackage{latexsym,amssymb,amsmath}
\usepackage{hyperref}
\usepackage{comment}
\usepackage{mathrsfs}
\usepackage{amsmath,amscd}
\usepackage[enableskew]{youngtab}
\usepackage[all]{xy}
\usepackage{empheq}

\usepackage{amsmath}

\usepackage{empheq}
\usepackage{xcolor}
\definecolor{lightgreen}{HTML}{90EE90}

\newcommand{\irrcc}{\irr(\cc)}

 \newcommand{\fpcc}{\fp(\cc)}

\newcommand{\bxt}{\boxtimes}

\newtheorem{theorem}{Theorem}[section]

\newtheorem{proposition}[theorem]{Proposition}
\newtheorem{corollary}[theorem]{Corollary}

\newcommand{\vect}{\mtr{Vec}}

\newcommand{\sent}{\mapsto}

\newcommand{\muj}{\mu_j}

\newcommand\C{\mathcal{C}}
\DeclareMathOperator{\id}{id}

\DeclareMathOperator{\rank}{{rank}}

\excludecomment{verlong}
\includecomment{vershort}
\excludecomment{noncompile}

\usepackage{tikz}
\usetikzlibrary{matrix}
\usepackage{combelow}
\newcommand{\ccb}{{\mathcal B}}

\usepackage{amsmath,amsthm,amsfonts,amssymb}
\usepackage{amsmath, amsthm, amssymb, amscd, amsfonts}
\usepackage[all]{xy}
\usepackage{amssymb, amsthm, amsmath}
\usepackage{amsfonts}
\usepackage{color}
\DeclareMathOperator{\ad}{ad} 

\newcommand{\ra}{\rightarrow}

\newcommand{\ot}{\otimes}

\newcommand{\xra}{\xrightarrow}
\newcommand{\mtc}{\mathcal}

\newcommand{\lam}{\lambda}

\newcommand{\Lam}{\Lambda}

\newcommand{\al}{\alpha}
\newcommand{\eps}{\epsilon}

\newcommand{\wte}{{\widehta{E}}}

\newcommand{\ul}{\underline}

\newcommand{\lh}{\leftharpoonup}
\newcommand{\whb}{{\widehat{(H, \mtc B)}}}
\newcommand{\hb}{{(H, \mtc B)}}
\newtheorem{thm}[theorem]{Theorem}
\newtheorem{rem}[theorem]{Remark}
\newtheorem{defn}[theorem]{Definition}
\newtheorem{lem}[theorem]{Lemma}
\theoremstyle{plain}

\newcommand{\ch}{\chi}
\newcommand{\mtr}{\mathrm}

\newcommand{\ncm}{\newcommand}
\ncm{\np}{\newpage}
\ncm{\ebl}{\end{thebibliography}}
\ncm{\bbl}{\begin{thebibliography}}
\ncm{\chd}{_{ _{\ch}}}
\ncm{\ald}{_{ _{\al}}}
\newcommand{\blam}{\Lam}
\ncm{\cP}{\mathcal{P}}
\ncm{\ei}{e_i}
\ncm{\eij}{e_{i,\;j}}
\ncm{\bt}{\begin{thm}}
\ncm{\bdef}{\begin{defn}}
\ncm{\edf}{\end{defn}}
\ncm{\et}{\end{thm}}
\ncm{\bc}{\begin{corollary}}
\ncm{\bl}{\begin{lem}}
\ncm{\el}{\end{lem}}
\ncm{\bpf}{\begin{proof}}
\ncm{\epf}{\end{proof}}
\ncm{\ec}{\end{corollary}}
\ncm{\ord}{\mtr{ord}}
\ncm{\er}{\end{rem}}
\ncm{\br}{\begin{rem}}
\ncm{\bn}{\begin}

\ncm{\bp}{\begin{proposition}}
\ncm{\ep}{\end{proposition}}
\ncm{\bd}{
\begin{document}}
\ncm{\ed}{\end{document}}
\ncm{\beq}{\begin{equation}}
\ncm{\beqn}{\begin{equation*}}
\ncm{\eeq}{\end{equation}}
\ncm{\eeqn}{\end{equation*}}
\ncm{\bea}{\begin{eqnarray}}
\ncm{\eea}{\end{eqnarray}}
\ncm{\beanon}{\begin{eqnarray*}}
\ncm{\eeanon}{\end{eqnarray*}}\ncm{\ek}{\eps|_K}\ncm{\diez}{\#}
\ncm{\bwt}{\bowtie}
\ncm{\cC}{\mtc{C}}\ncm{\cc}{\mtc{C}}
\ncm{\cX}{\mtc{X}}
\ncm{\wt}{\widetilde}
\ncm{\sg}{\sigma}
\ncm{\Rep}{\mathrm{Rep}}
\ncm{\Aut}{\mathrm{Aut}}
\DeclareMathOperator{\Irr}{Irr}
\ncm{\X}{\mathcal{X}}
\ncm{\cA}{\mathcal{A}}
\ncm{\HKer}{\mtr{HKer}}
\ncm{\LKER}{\mtr{LKer}}
\ncm{\aad}{\mtr{ad}}
\newcommand{\mbf}{\mathbb F}
\ncm{\Dr}{\mtr{D}}
\ncm{\cD}{{\mathcal{D}}}\ncm{\cd}{{\mathcal{D}}}\ncm{\ce}{{\mathcal{E}}}
\ncm{\G}{\mathcal{G}}
\ncm{\Dc}{\mtc{D}}
\ncm{\E}{\mtc{E}}
\ncm{\fp}{\mtr{FPdim}}
\ncm{\Vc}{\mtr{Vec}}
\ncm{\cK}{\mtc{K}}
\ncm{\cM}{\mtc{M}}
\ncm{\cE}{\mtc{E}}
\ncm{\cS}{\mtc{S}}

\newcommand{{\ipr}}{i'}
\newcommand{\tomega}{\widetilde{\omega}}

\DeclareMathOperator{\End}{End}
\ncm{\cop}{\mtr{cop}}
\ncm{\op}{\mtr{op}}
\ncm{\chr}{character }\ncm{\ck}{\mtc{K}}
\ncm{\bw}{\bwt}
\ncm{\hker}{\mtr{HKer}}
\ncm{\bx}{\boxtimes}
\ncm{\blue}{\textcolor[rgb]{.00, .00, 1.00}}
\ncm{\bluer}{\textcolor[rgb]{.30, .30, .70}}
\ncm{\red}{\textcolor[rgb]{1.00, .00, .00}}
\ncm{\green}{\textcolor[rgb]{.50, 0.20, .90}}
\ncm{\bne}{\begin{enumerate}}
\ncm{\ene}{\end{enumerate}}
\ncm{\lker}{\mtr{LKer}}
\ncm{\md}{\medbreak}
\ncm{\rep}{\Rep}\ncm{\ind}{\mtr{ind}}
\ncm{\mdn}{\md\noindent}
\ncm{\dd}{$}
\ncm{\up}{^}
\newcommand{\tcs}{\text}
\newcommand{\mbb}{\mathbb B}
\newcommand{\vs}{\mathbb V}
\newcommand{\sth}{suppose that\;}
\newcommand\rad{\operatorname{rad}}
\newcommand{\itm}{\item}
\newcommand{\dbd}{$$}
\newcommand{\mol}{\mtr{mod}}
 \newcommand{\ro}{\rho}
\newcommand{\irr}{\mathrm{Irr}}
\newcommand{\mbc}{\mathbb C}
\newcommand{\mbs}{\mathbb S}
\newcommand{\mbz}{\mathbb Z}
\newcommand{\ct}{\mtc T}
\newcommand{\sm}{\setminus}
\newcommand{\epl}{^{+}}
\newcommand{\sbsq}{\subseteq}
\newcommand{\sbs}{\subset}
\newcommand{\cco}{\mtr{co}}
\newcommand{\cz}{\mathcal{Z}}
\newcommand{\dual}{^{*}}
\newcommand{\Gm}{\Gamma}
\ncm{\cY}{\mtc{Y}}
\newcommand\ZZ{{\mathbb Z}} 
\newcommand{\bab}{\color{DarkOrchid}{}}
\newcommand{\eab}{\normalcolor{}}
\newcommand{\subs}{\subsection}
\newcommand{\cv}{\mtc{V}}
  \newcommand{\grn}{\green}
\newcommand{\dt}{\delta}

\newcommand{\ccf}{\mathrm{ {CF}(\cc)}}
\newcommand{\cce}{\mathrm{ {CE}(\cc)}}
\newcommand{\cecc}{\mathrm{ {CE}(\cc)}}
\newcommand{\cecd}{\mathrm{ {CE}(\cd)}}
\newcommand{\kk}{\Bbbk}
\newcommand{\otL}{\ot_{L}}
\newcommand{\otl}{\ot_{L}}
\newcommand{\unpsi}{1_{\psi}}
\newcommand{\epsi}{e_{\psi}}
\newcommand{\ephi}{e_{\phi}}
\newcommand{\ech}{e_{\ch}}
\newcommand{\nleftcid}{\text{left normal  coideal subalgebra}}
\newcommand{\dimL}{\dim_{\kk}L}
\newcommand{\cl}{\mtc L}
\newcommand{\mj}{\mtc J}
\newcommand{\tl}{\tilde L}
\newcommand{\tL}{\tilde L}
\newcommand{\tpsi}{\tilde(\psi)}
\newcommand{\tmx}{\tilde{\mtc X}}
\newcommand{\zlh}{\mathrm{ZL}}
\newcommand{\ba}{\mathrm A}
\newcommand{\bv}{\mathrm V}
\newcommand{\zhopf}{\mtc{Z}_{\mtr{Hopf}}}
\newcommand{\lstar}{L^{*}}
\newcommand{\ldstar}{L^{**}}
\newcommand{\mstar}{M^{*}}
\newcommand{\mdstar}{M^{**}}
\newcommand{\lkera}{\lker_{A}}
\newcommand{\mdprime}{M''}
\newcommand{\ldprime}{L''}
\newcommand{\cm}{\mtc M}
\newcommand{\ccm}{\mathcal M}
\newcommand{\cn}{\mathcal N}
\newcommand{\ccn}{\mathcal N}
\newcommand{\rx}{\mtr{Rex}}
\newcommand{\cca}{\ca}
\newcommand{\ih}{\underline{\mtr{Hom}}}
\newcommand{\cih}{\underline{\mtr{coHom}}}
\newcommand{\hm}{\mtr{ {Hom}}}
\newcommand{\cov}{\mtr{coev}}
\newcommand{\rora}{\rho^{\mtr{ra}}}
\newcommand{\rola}{\rho^{\mtr{la}}}
\newcommand{\cx}{\mtc X}
 \newcommand{\cZ}{\cz}
 \newcommand{\ca}{\cA}
 \newcommand{\stat}{\noindent}
 \newcommand{\bfa}{{\bf A}}
 \newcommand{\unu}{\mathbf{1}}
 \newcommand{\barzu}{{\bar {  Z}(\unu)}}
 
\newcommand{\idx}{\id_X}
\newcommand{\lprime}{L'}
\newcommand{\mprime}{M'}
\newcommand{\nat}{ \mtr{{  Nat}}}
\newcommand{\ft}{\mtc F_\lam}
\newcommand{\rhau}{\rightharpoonup}
\newcommand{\lhau}{\leftharpoonup}
\newcommand{\cf}{\mathrm{ {CF}}}

\newcommand{\cfc}{\mathrm{{CF}}(\cc)}
\newcommand{\csu}{\overline{\mathfrak{  C}}}
\newcommand{\cfcc}{{\mathrm{CF}(\cc)}}
\newcommand{\catfcc}{\mathrm{ {CF}}(\cc)}
\newcommand{\cfcd}{\mathrm{CF}(\cd)}
\newcommand{\cfd}{\mathrm{CF}(\cd)}
\newcommand{\czcc}{{\cz(\cc)}}
\newcommand{\czcd}{{\cz(\cd)}}
\newcommand{\czt}{{\cz(\cz(\cc))}}
\newcommand{\enx}{\mtr{  End}}
\newcommand{\runu}{R(\unu)}

\newcommand{\bdfn}{\bn{defn}}
\newcommand{\edfn}{\end{defn}}
\newcommand{\deltax}{\delta_X}
\newcommand{\deltav}{\delta_V}
\newcommand{\repcca}{\rep_\cc(A)}
\newcommand{\xotay}{X \ot_A Y}
\newcommand{\xoty}{X \ot Y}
\newcommand{\votw}{V \ot W}
\newcommand{\votaw}{V \ot_A W}
\newcommand{\dimax}{\dim_AX}
\newcommand{\dimccx}{\dim_\cc(X)}
\newcommand{\dimcca}{\dim_\cc(A)}
\newcommand{\dimccv}{\dim_\cc(V)}
\newcommand{\dima}{\dim_A}
\newcommand{\biga}{A}
\newcommand{\comp}{\mathbb C}
\newcommand{\tehtaa}{\theta_A}
\newcommand{\tetaa}{\theta_A}
\newcommand{\ida}{\id_A}
\newcommand{\hma}{\hm_A}
\newcommand{\hmcc}{\hm_\cc}
\newcommand{\fv}{F(V)}
\newcommand{\fw}{F(W)}
\newcommand{\ota}{\ot_A}
\newcommand{\repza}{\rep_\cc^0(A)}
\newcommand{\epsa}{\eps_A}
\newcommand{\bndefn}{\bn{defn}}
\newcommand{\edefn}{\end{defn}}
\newcommand{\bdefn}{\bn{defn}}

\newcommand{\vld}{V^{*}}
\newcommand{\vldd}{V^{**}}
\newcommand{\xld}{X^{*}}
\newcommand{\xldd}{X^{**}}
\newcommand{\yld}{Y^{*}}
\newcommand{\yldd}{Y^{**}}
\newcommand{\aldu}{A^{*}}
\newcommand{\aldd}{A^{**}}

\newcommand{\ia}{\mtr{i}_A}
\newcommand{\aota}{A\ot A}

\newcommand{\idv}{\id_V}

\newcommand{\ld}{^*}
\newcommand{\repg}{\rep(G)}

\newcommand{\thetav}{\theta_V}

\newcommand{\tta}{\theta_A}

\newcommand{\muv}{\mu_V}
\newcommand{\muw}{\mu_W}

\newcommand{\dimcc}{\dim(\cc)}
\newcommand{\chii}{\chi_i}
\newcommand{\chistar}{\ch_{i^*}}
\newcommand{\chj}{\ch_j}
\newcommand{\chm}{\ch_m}
\newcommand{\chn}{\ch_n}
\newcommand{\dimvi}{\dim(V_i)}
\newcommand{\mtcd}{Q}
\newcommand{\mtca}{\mtc A}
\newcommand{\lamcd}{\lam_\cd}
\newcommand{\fpdimcd}{\fp(\cd)}
\newcommand{\laml}{\lam_L}
\newcommand{\apm}{A//M}
\newcommand{\apl}{A//L}
\newcommand{\repapm}{\rep(\apm)}
\newcommand{\repapl}{\rep(\apl)}
\newcommand{\dimvj}{\dim(V_j)}
\newcommand{\dvi}{\dim(V_i)}
\newcommand{\dvj}{\dim(V_j)}
\newcommand{\sumjtom}{\sum_{j=0}^m}
\newcommand{\sumitom}{\sum_{i=0}^m}
\newcommand{\sij}{s_{ij}}
\newcommand{\sji}{s_{ji}}
\newcommand{\dxj}{d_j}
\newcommand{\dxi}{\di }
\newcommand{\dimka}{\dim_{\kk}(A)}
\newcommand{\dimk}{\dim_{\kk}}
\newcommand{\blaml}{\blam_L}
\newcommand{\sumjtor}{\sum_{j=0}^r}
\newcommand{\dimkl}{\dim_{\kk}(L)}
\newcommand{\mtcjl}{\mtc J_L}
\newcommand{\vota}{ V\ot A}
\newcommand{\vi}{V_i}
\newcommand{\vj}{V_j}
\newcommand{\dimcd}{\dim(\cd)}

\newcommand{\alij}{{\al_{ _{ij}}}}
\newcommand{\alji}{{\al_{ _{ji}}}}
\newcommand{\rcc}{r_\cc}
\newcommand{\rcd}{r_\cd}
\newcommand{\clsx}{[X]}
\newcommand{\clsy}{[Y]}
\newcommand{\clsz}{[Z]}
\newcommand{\rcdp}{r_{\cd'}}
\newcommand{\sumjtorp}{\sum_{j=0}^{r'}}
\newcommand{\aljm}{{\al_{ _{jm}}}}
\newcommand{\aljn}{{\al_{ _{jn}}}}
\newcommand{\sjm}{s_{jm}}
\newcommand{\smj}{s_{mj}}
\newcommand{\snj}{s_{nj}}

\newcommand{\betaij}{\beta_{ _{ij}}}
\newcommand{\betaji}{\beta_{ _{ji}}}
\newcommand{\gammaij}{\gamma_{ _{ij}}}
\newcommand{\gammaji}{\gamma_{ _{ji}}}
 \newcommand{\ip}{i'}
\newcommand{\sumjtoprp}{\sum_{j=0}^{r'}}
\newcommand{\sumjtopr}{\sum_{j=0}^{r}}
 \newcommand{\teh}{\tilde{h}}
\newcommand{\cdp}{{\cd'}}\newcommand{\xphii}{X_{\phi(i)}}
\newcommand{\inv}{^{-1}}

\newcommand{\fq}{\mtr f_{ Q}}
\newcommand{\tr}{\mtr{tr}}
\newcommand{\rtwone}{R_{21}R}

\newcommand{\ccad}{{\cc_{\mtr{ad}}}}
\newcommand{\ccpt}{{\cc_{\mtr{pt}}}}
\newcommand{\qtr}{quasi-triangular\;}
\newcommand{\trq}{\tr_q}

\newcommand{\repal}{\mtr{Rep}(A//L)}
\newcommand{\lkeravi}{\lker_A(V_i)}
\newcommand{\lkeravj}{\lker_A(V_j)}
\newcommand{\cross}[1][1pt]{\ooalign{%
 \rule[1ex]{1ex}{#1}\cr
 \hss\rule{#1}{.7em}\hss\cr}}
\newcommand{\blml}{\blam_L} 
\newcommand{\phir}{\phi_R}
\newcommand{\kda}{{  \Phi(A)}}

\newcommand{\mtcil}{\mtc{I}_L}

\newcommand{\un}{\unu}
\newcommand{\tfl}{\mtc{T}}
\newcommand{\barzm}{\barz(M)}
\newcommand{\barzn}{\barz(N)}
\newcommand{\ccr}{\mtc R^{\cc}}
\newcommand{\ulc}{\ul{\cc}}

\newcommand{\pimx}{\pi_{M;\;X}}
\newcommand{\pinx}{\pi_{N;\;X}}
\newcommand{\acc}{{\mathrm A_\cc}}
\newcommand{\epsu}{\eps_\unu}

\newcommand{\ob}{\mtr{Obj}}
\newcommand{\obc}{\mtr{Obj(\cc)}}
\newcommand{\ccop}{\cc^{\mtr{op}}}
\newcommand{\mtf}{\mtc F_\lam}
\newcommand{\mtfi}{\mtc F^{-1}_\lam}
\newcommand{\elcd}{\ell_\cd}
\newcommand{\mcid}{\mtc I_\cd}
\newcommand{\mcidp}{\mtc I_{\cd'}}
\newcommand{\wtildelcd}{\widetilde{\elcd}}
\newcommand{\wtildelcdp}{\widetilde{\ell_{\cd'}}}
\newcommand{\cpt}{\cc_{\mtr{pt}}}
\newcommand{\barzr}{\barz_\cd}
\newcommand{\barzv}{\barz(V)}
\newcommand{\acd}{\mathrm A_\cd}
\newcommand{\czrcd}{\cz_\cc(\cd)}
\newcommand{\sml}{\Small}
\newcommand{\bs}{\blue{\Small }}
\newcommand{\yd}{Yetter-Drinfeld}

\newcommand{\sumitor}{\sum_{i=0}^r}
\newcommand{\cdop}{\cd^{\mtr{op}}}
\newcommand{\ccrev}{\cc^{\mtr{rev}}}
\newcommand{\barz}{{\bar{\mathrm Z}}}
\newcommand{\etl}{etale\;}
\newcommand{\czca}{\cz(\ca)}

\newcommand{\tetx}{\text}
\newcommand{\widehta}{\widehat}
\newcommand{\wdhat}{\widehat}
\newcommand{\wht}{\widehat}
\newcommand{\cofa}{{\mathbb C[\mtc B]}}
\newcommand{\wdt}{\widehat}
\newcommand{\dl}{{^\#}}
\newcommand{\comx}{\mathbb C}

\newcommand{\sgj}{{\sg(j)}}

\newcommand{\mujo}{\mu_\jo}
\newcommand{\mujtw}{\mu_\jtw}
\newcommand{\adz}{a^{\#}}
\newcommand{\bdz}{b^{\#}}

\newcommand{\spr}{S^\perp}
\newcommand{\cofs}{\comp [S]}
\newcommand{\spz}{S^{\perp_z}}

\newcommand{\omz}{\omega_z}
\newcommand{\zg}{\mathrm{Z}(S)}
\newcommand{\aling}{{\al \in g}}

\newcommand{\blkg}{\mtr{Bl}(g)}
\newcommand{\clsg}{\mtr{Cl}(g)}
\newcommand{\mtadinv}{\mtc G^{{-1}}}
\newcommand{\muk}{\mu_{k}}
\newcommand{\mta}{\mtc F}
\newcommand{\cofad}{\comp[\wdht A]}
\newcommand{\wtau}{\wdht{\tau}}
\newcommand{\mtainv}{{\mta}^{-1}}
\newcommand{\wdht}{\widehat}
\newcommand{\augm}{\mtr{aug}}
\newcommand{\mua}{\wdht {\wdht a}}
\newcommand{\aps}{A//S}
\newcommand{\cssa}{\cc(S, A)}
\newcommand{\aug}{\mtr{aug}}
\newcommand{\rss}{{\big|_S}}
\newcommand{\gprp}{g^\perp}
\newcommand{\alins}{{s \in S}}

\newcommand{\sz}{s^{D}}
\newcommand{\wmu}{\widehta{\mu}}
\newcommand{\wmui}{\widehta{\mu}_i}
\newcommand{\wmuj}{\widehta{\mu}_j}
\newcommand{\wch}{\widehta{\ch}}

\newcommand{\wzd}{\widehat{d}}
\newcommand{\wpm}{\widehat{P}}
\newcommand{\wps}{\widehat{p}}

\newcommand{\gal}{\mtr{Gal}}
\newcommand{\galkq}{\gal(\mathbb K/\mathbb Q)}
\newcommand{\sgh}{\sg_{ _{H}}}
\newcommand{\sggi}{{\sg(i)}}
\newcommand{\sge}{\sg_{_{\widehat R}}}
\newcommand{\unue}{{\unu_{\cecc}}}

\newcommand{\mtcf}{\mtc {F}}

\newcommand{\wsgf}{\widehat{{\sg}_{ _F}}}
\newcommand{\sghstar}{{{\sg}_{ _{H^*}}}}
\newcommand{\we}{\widehta{E}}
\newcommand{\sumktom}{\sum_{k=0}^m}

\newcommand{\wf}{\widehat{F}}

\newcommand{\hsgj}{\widehat{\sg}(j)}
\newcommand{\whsgi}{\widehta{\sg}(i)}

\newcommand{\wpp}{\widehat{p}}
\newcommand{\tauj}{{\tau(j)}}
\newcommand{\dimcctauj}{\dim(\cc^\tauj)}
\newcommand{\etas}{{\eta(s)}}
\newcommand{\mcc}{m_H}

\newcommand{\wal}{\widehta{\al}}
\newcommand{\wj}{\widehat{\mtc J}}
\newcommand{\galc}{\mtr{Gal}_{\cc}}
\newcommand{\galz}{\mtr{Gal}_{\czcc}}
\newcommand{\wjr}{\widehat{J}_{R}}

\newcommand{\dimcck}{\dim(\cc^k)}

\newcommand{\wgrcc}{\widehat{\mtr{Gr}(\cc)}}
\newcommand{\nchi}{{\frac{\ch_i}{\di }}} \newcommand{\nchj}{{\frac{\ch_j}{\dxj}}}
\newcommand{\wni}{{\widehat{n}_i}}

\newcommand{\sgte}{\widetilde{\sg_E}}

\newcommand{\mtad}{\mtc G}
\newcommand{\whj}{\widehta{h}_j}
\newcommand{\jdl}{{j\dl}}
\newcommand{\wcfcc}{\widehat{\cfcc}}
\newcommand{\mutauj}{\mu_{\tau(j)}}
\newcommand{\tauk}{\tau(k)}
\newcommand{\muzm}{{\mu_0^{-}}}
\newcommand{\sqrtog}{\sqrt{|G|}}
\newcommand{\muz}{\mu_0}
\newcommand{\njtw}{n_\jtw}
\newcommand{\njo}{n_\jo}
\newcommand{\fjo}{F_\jo}
\newcommand{\fjtw}{F_\jtw}
\newcommand{\wta}{\widehat{A}}

\newcommand{\dol}{{^{\circ}}}
\newcommand{\bdl}{{b\dl}}
\newcommand{\jdol}{{j\dol}}
\newcommand{\fj}{F_j}

\newcommand{\cwta}{\comp[\wta]}

\newcommand{\hx}{\widehta{x}}
\newcommand{\hy}{\widehta{y}}

\newcommand{\cal}{\mtc A_{\al}}
\newcommand{\xuu}{x_{uu}}
\newcommand{\wxuu}{\widehat{\xuu}}
\newcommand{\xvv}{x_{vv}}
\newcommand{\xuv}{x_{uv}}
\newcommand{\xmn}{x_{m,n}}
\newcommand{\buvmn}{B^{u,v}_{m,n}}
\newcommand{\blm}{\blam}
\newcommand{\dimccr}{\dim(\cc^r)}
\newcommand{\adl}{a\dl}
\newcommand{\sumltom}{\sum_{l=0}^m}

\newcommand{\mbq}{\mathbb Q}
\newcommand{\mbqs}{\mathbb Q(S)}
\newcommand{\mbk}{\mathbb K}
\newcommand{\mz}{\mathbb Z}

\newcommand{\wsgj}{\widehat{\sigma}(j)}
\newcommand{\wsgi}{\widehat{\sigma}(i)}
\newcommand{\wg}{\widehat{g}}
\newcommand{\wtf}{\widehat{F}}
\newcommand{\galqspq}{\mtr{Gal}(\mathbb Q(S)/\mathbb Q)}
\newcommand{\cctauj}{\cc^{\tau(j)}}
\newcommand{\cctauk}{\cc^{\tau(k)}}
\newcommand{\wtfj}{\widetilde{F_j}}
\newcommand{\wfj}{\widetilde{F_j}}
\newcommand{\wtmuj}{\widetilde{\mu_j}}
\newcommand{\wmtcfj}{{\widetilde{\mtc F}_j}}
\newcommand{\mtfr}{\mtr{F_a}}
\newcommand{\wdr}{R_\comp^*}

\newcommand{\fgph}{{F_{G/H}}}
\newcommand{\wcfj}{\wmtcfj}

\newcommand{\nxi}{{\frac{x_i}{\di }}}
\newcommand{\fpr}{{\fp(R)}}
\newcommand{\nxs}{{\frac{x_s}{d_s}}}

 \newcommand{\mtfme}{\mtc F}
\newcommand{\chic}{\ch_i^{\circ}}
\newcommand{\chjc}{\ch_j^{\circ}}
\newcommand{\mtfsh}{{\mtc F_\lam}}
\newcommand{\mupq}{{\mu_{pq}}}

 \newcommand{\tlam}{{\widetilde{\lam}}}
 \newcommand{\chid}{{\ch_i^{\circ}}}
\newcommand{\rc}{{R_\comp}}
\newcommand{\rgo}{{\mathbb R_{\geq 0}}}
\newcommand{\sumrorc}{{{\sum\limits_{\ro \in \rc}}}}
\newcommand{\aliro}{{\ro(x_i)}}

\newcommand{\barjd}{{\bar{\mtc J_\cd}}}
\newcommand{\lbarcj}{{\frac{C_j}{{\dim(\mathcal C^j)}}}}
\newcommand{\omtcb}{{\overline{\mathcal B}}}
\newcommand{\whr}{{\widehat{R}}}
\newcommand{\nxj}{\frac{x_j}{\dxj}}
\newcommand{\nxk}{\frac{x_k}{d_k}}
\newcommand{\onkij}{{\overline{N^k_{ij}}}}
\newcommand{\sgk}{{\sigma(k)}}
\newcommand{\sgl}{{\sigma(l)}}
\newcommand{\fqi}{{\fq^{-1}}}
\newcommand{\wdb}{{\widehat{\mtc B}}}
\newcommand{\mtcb}{{\mtc B}}

\newcommand{\nif}{{h_i}}
\newcommand{\rb}{(R, \mtc B)}
 \newcommand{\nxip}{{\frac{x_{i'}}{d_{i'}}}}
\newcommand{\mujp}{{\mu_{j'}}}

\newcommand{\etai}{{\eta(i)}}
\newcommand{\wsg}{{\widehat{\sg}}}
\newcommand{\wsgh}{{\wsg_{ _{H}}}}
\newcommand{\wtaui}{{\widehat{\tau}(i)}}
\newcommand{\wsghstar}{{\wsg_{H^*}}}
\newcommand{\wtauj}{{\wtau(j)}}
\newcommand{\weta}{{\widehat{\eta}}}
\newcommand{\detai}{{d_{ _{\eta(i)}}}}

\newcommand{\hbz}{{(H, \mtc B, \mu_0)}}
\newcommand{\whbz}{{\widehat{\hbz}}}
\newcommand{\tsgh}{{{\widetilde{\sgh}}}}
\newcommand{\ghb}{{G\hb}}
\newcommand{\taujo}{{\tau(\jo)}}
\newcommand{\taujtw}{{\tau(\jtw)}}
\newcommand{\mutauk}{{\mu_{\tauk}}}

\newcommand{\hetai}{{h_{ _\etai}}}
\newcommand{\xetai}{{x_{ _\etai}}}
\newcommand{\wn}{{\widehat{n}}}
\newcommand{\wh}{{\widehat{h}}}
\newcommand{\distar}{{d_{i^*}}}
\newcommand{\dwtaui}{{d_{ _{\wtaui}}}}
\newcommand{\sumiptom}{{\sum_{\ip=0}^m}}
\newcommand{\alitaugj}{{\al_{ _{i\tau_g(j)}}}}
\newcommand{\dtaugj}{{d_{ _{\tau_g(j)}}}}
\newcommand{\mip}{{M(\ip)}}
 \newcommand{\sumttom}{{\sum_{t=0}^m}}
 \newcommand{\muxi}{{\mu_{ _{[X_i]}}}}
\newcommand{\muxip}{{\mu_{ _{[X_{\ip}]}}}}
\newcommand{\ncj}{{\frac{C_j}{\dim(\cc^j)}}}
\newcommand{\jp}{{j'}}
\newcommand{\minv}{{M^{-1}}}
\newcommand{\tfq}{{\widehat{\fq}}}
\newcommand{\catcecc}{{\mtr{CE}(\cc)}}
\newcommand{\wcatfcc}{{\widehat{\catfcc}}}
\newcommand{\mforall}{{\;\;\text{for all}\;\;}}
\newcommand{\what}{\widehat}
\newcommand{\wfz}{{\widehat{F}_0}}
\newcommand{\hbfr}{{(H, \mtc B, \fp)}}
\newcommand{\sgn}{{\mtr{sgn}}}
\newcommand{\wir}{{\widehat{R}}}
\newcommand{\gcc}{{G(\cc)}}
\newcommand{\jccpt}{{J_{ _{\ccpt}}}}
\newcommand{\jccad}{{J_{ _{\ccad}}}}
\newcommand{\fpccad}{{\fp(\ccad)}}
\newcommand{\fpccpt}{{\fp(\ccpt)}}
\newcommand{\kc}{{K(\cc)}}
\newcommand{\wkc}{{\widehat{\kc}}}

\newcommand{\hs}{{(L,\;\mtc S)}}
\newcommand{\rbad}{{H_{ _{ad}}}}
\newcommand{\hbad}{{\hb_{ad}}}
\newcommand{\jhbad}{{J_{\hbad}}}
\newcommand{\htt}{{(K, \mtc T)}}
\newcommand{\jhtt}{{\mtc J_{ _{\htt}}}}
\newcommand{\jhs}{{\mtc J_{ _{\hs}}}}
\newcommand{\lamhs}{{\lam_{ _{\hs}}}}
\newcommand{\lamhtt}{{\lam_{ _{\htt}}}}

\newcommand{\coo}{{co}}
\newcommand{\wdhad}{{(\wdh)_{ad}}}
\newcommand{\had}{{H_{ad}}}
\newcommand{\wdh}{\widehat{H}}
\newcommand{\whbad}{{\whb_{ _{ad}}}}
\newcommand{\kerhb}{{\ker_{ _{\hb}}}}
\newcommand{\gwdh}{{G(\wdh)}}
\newcommand{\nxl}{\frac{x_l}{d_l}}

\newcommand{\proditom}{{\prod_{i=0}^m}}

\newcommand{\wdhn}{{\wdh^{(n)}}}
\newcommand{\hn}{{H_{(n)}}}
\newcommand{\nox}{{\frac{x}{\fp(x)}}}
\newcommand{\mujstar}{{\mu_{j^\#}}}
\newcommand{\sco}{{S^\coo}}
\newcommand{\rrad}{{I(1)}}
\newcommand{\istar}{{i^*}}
\newcommand{\mtcs}{{\mtc S}}

\newcommand{\qghb}{{\ghb}}
\newcommand{\wqghb}{{{G\whb}}}
\newcommand{\nxp}{{\frac{x_p}{\fp(x_p)}}}
\newcommand{\nxm}{{\frac{x_m}{\fp(x_m)}}}

\newcommand{\whp}{{\widehat{P}}}
\newcommand{\prodjtom}{{\prod_{j=0}^m}}
\newtheorem{statement}[theorem]{Statement}
\newcommand{\wghb}{{G\whb}}
\newcommand{\gwhb}{{G\whb}}
\title
{On Harada's identity  and some other consequences of Burnside's vanishing property}

\author{Sebastian Burciu}
\address{Inst.\ of Math.\ ``Simion Stoilow" of the Romanian Academy P.O. Box 1-764, RO-014700, Bucharest, Romania}
\email{sebastian.burciu@imar.ro}
\newcommand{\fpcctw}{{\fpcc_{\geq 2}}}
\newcommand{\fpccone}{{\fpcc_{1}}}
\newcommand{\wkcad}{{\widehat{\kc}_{ad}}}
\date{\today}
\subjclass{18M20; 16T30}
\bd
\thanks{The first author is supported by a grant of the Ministry of Research, Innovation and Digitization, CNCS/CCCDI - UEFISCDI, project number PN-III-P4-ID-PCE-2020-0878, within PNCDI III}
\maketitle
\begin{abstract}
In this short note, we prove some consequences of Burnside's vanishing property \cite{b-pa}. It is known that Harada's identity concerning the product of all conjugacy classes of a finite group is a consequence of Burnside's vanishing property of characters.  We prove a similar formula for any weakly-integral fusion category. In the second part, we prove some structural results concerning nilpotent and modular fusion categories.
As an application we show that any fusion category of dimension $p^2q^2r^2d$ with $p<q<r$ prime numbers and $d$ a square-free integer is weakly-group theoretical. 
\end{abstract}

\setcounter{tocdepth}{1}
\tableofcontents
\section{Introduction}\label{introd}
A classical vanishing result of Burnside states that any  irreducible non-linear character  of a finite group vanishes on at least one element of the group. 

Based on this vanishing result, in \cite[Theorem 1]{harada}, K. Harada proved that the product of all conjugacy classes of a finite group form a coset with respect to the commutator subgroup $G'$. 

\bt (Harada, \cite{harada}) For any finite group $G$, the following identity holds:
$$
\left(\prod_{j=0}^m \frac{C_j}{|\mathcal C^j|}\right)^2=\frac{1}{|G'|}\sum_{\mathcal C^j\subseteq G'} C_j
$$
where $G'$ is the commutator subgroup of $G$. 
\et
In the above theorem $\{\mathcal C^j\}_{j=0}^m$ are all the conjugacy classes of $G$ and $|\mathcal C^j|$ is the size of the conjugacy class. The elements $C_j:=\sum_{g \in \mathcal C^j}g$ are the class sums associated to $\mathcal C^j$.

One of the goals of this short note is to present a similar formula for any weakly-integral fusion category.  For weakly-integral fusion categories Burnside's vanishing result was extended in \cite{b-galois}. Some other consequences of Brauer's vanishing result at the level of fusion categories are also presented.

Let $\cc$ be a weakly-integral fusion category, $\cc^j$ the set of all its conjugacy classes and $C_j$ their associated class sums. Our first main result is the following:
\bt\label{harada}
For any weakly-integral fusion category, with the above notations,  one has 
$$
(\prod_{j=0}^m \frac{C_j}{\fp(\cc^j)})^2=\frac{\fp(\ccpt)}{\fpcc}\big(\sum_{\muj\in J_{\ccpt}}C_j\big)
$$
\et
For the definition of conjugacy classes $\cc^j$ and their associated class sums $C_j$, see Section \ref{prelim}. Moreover, the support of the fusion subcategory $\ccpt$ is defined in \cite[Subsection 4.2]{ccc-march}.


 In \cite[Theorem 4.2]{o-yu} it was shown that any weakly group-theoretical modular tensor category whose dimension is $nd$ with $(n,d)=1$ is in fact a Deligne product of $\cc(\mathbb Z_q, d) \boxtimes \cd$ of two other modular tensor categories.

Suppose that $\cc$ is a weakly-integral modular category. As usually we suppose that $\cc$ is of type $(d_1,n_1;d_2,n_2;\dots d_r,n_r)$ with $d_1$=1. Our second main result of this note is the following:
\bt\label{div}
Suppose that $\cc$ is a weakly-integral modular tensor category with $\fpcc=nd$ with $(n,d)=1$ and $d$ a square-free integer. 

Then, with the above notations one has  $d\big| n_i$ for all $i\geq 1$.
\et
As an application we prove the following:
\bt\label{pqr}
Let $\cc$ be an integral modular fusion category of dimension $p^2q^2r^2d$ with $d$ a square-free number and $p<q<r$ three prime numbers. Then any such fusion category is weakly group-theoretical and can be written as a Deligne product $\cca\bxt \ccb$ with $\cca\simeq \cc(\mathbb Z_q, d)$ and $\ccb$ a weakly group theoretical fusion category of dimension $p^2q^2r^2$.
\et
Note that the case $d=1$ of the above result was treated in \cite[Prop. 3.8.]{j-plav-odd}. Proposition \ref{dim:nilp} proves another consequence of  Burnside vanishing result for nilpotent fusion categories.

Shortly, this note is organized as follows:
In Section \ref{prelim} we remind the basics on pivotal fusion categories that are needed through the paper. In Section \ref{pf-harada} we prove Harada's identity for weakly-integral fusion categories. In Section \ref{pf-mtc} we prove Theorem \ref{div} and Proposition \ref{dim:nilp}
\section{Preliminaries}\label{prelim}

Throughout this paper we work over the ground field $\mathbb C$. For the basic theory of fusion categories, we refer the reader to the monograph \cite{EGNO15}. 

Recall that a pivotal structure of a rigid monoidal category $\cc$ is an isomorphism $j:\id_\cc\ra ()^{**}$ of monoidal functors. A pivotal monoidal category is a rigid monoidal category endowed with a pivotal structure.

For any fusion category $\cc$, the monoidal center (or the Drinfeld center) of $\cc$ is a braided fusion category $\czcc$ endowed with a forgetful functor $F:\czcc\ra \cc$, see \cite{EGNO15}.

The forgetful functor $F$ also admits a right adjoint functor $R:\cc \ra \czcc$  and $Z :=FR:\cc \ra \cc$ is a Hopf comonad defined as an end.  Indeed, following \cite[Section 2.6]{scalg} one has that 
\beq
Z(V)\simeq \int_{X\in \cc}X\ot V\ot X^*
\eeq
It is known that $A:=Z(\unu)$ has the structure of central commutative algebra in $\cz(\cc)$.

The vector space $\cecc:= \hm_{\C}(\unu, A) $ is called {\it the set of central elements.} There is a canonical bijection:
{\Small
\beq\label{canisom}
\cce\xra{\psi} \enx(\id_{\C}),\;\;\psi(a)_X=\ro_X(a\ot \\id_X),\;\; X=\unu\ot X\xra{a\ot \\id_X}X\ot Z(\unu)\xra{\ro_X} X
\eeq
}
The vector space $\cfcc:=\hm_\cc(A_\cc, \unu)$ is called the {\it space of class functions} of $\cc$. To any object $X\in \irr(\cc)$ Shimizu has attached an element $\ch_X\in \cfcc$ called {\it the character of $X$.}

Recall $R:\cc \ra \czcc$ is a right adjoint to the forgetful functor $F:\czcc \ra \cc$. As explained in \cite[Theorem 3.8]{scalg} this adjunction  gives an isomorphism of of $\comp$-algebras
\beq\label{adjisom}
\cfcc \xra{\cong} \mtr{End}_{\czcc}(R(\unu)),\;\; \ch\mapsto Z (\ch)\circ \delta_\unu.
\eeq
where $\delta_\unu$ is the comultiplication of the comonad $Z$ associated at $\unu$.

Let  $\lambda\in \cfcc$ be  a non-zero cointegral of $\cc$. The {\it Fourier transform} of $\cc$  associated to $\lambda$ is the linear map
\beq
\mtc F_{\lambda}:\cecc\ra \cfcc\;\;\text{given by}\;\;a \mapsto \lambda \lh \mtc S(a)
\eeq
\noindent
where $\mtc S: \cecc\ra \cecc$ is the above antipodal operator. The Fourier transform is a bijective $\kk$-linear map whose inverse is given in  \cite[Equation 5.17]{scalg}. Since the monoidal center $\czcc$ is also a fusion category we can write $R(\unu)=\bigoplus_{i=0}^m\mathcal C_i$ as a direct sum of simple objects in $\czcc$. Thus  $\mathcal C^{0},\dots, \mathcal C^{m}$ are the conjugacy classes of $\C$.  Since the unit object $\unu_{\czcc }$ is always a subobject of $R(\unu)$, we can assume $\mathcal C^{0} = \unu_{\czcc }$.

Let $\tilde F_0, \tilde F_1, \dots, \tilde F_m\in \enx_{\czcc}(R(\unu))$ be the canonical projections on each of these conjugacy classes. Let also $F_0, F_1, \dots F_m$ be also the corresponding primitive idempotents of $\cfcc$ under the canonical adjunction isomorphism $\cfcc \simeq \enx_{\czcc}(R(\unu))$ from Equation \eqref{canisom}.

We define $C_i:={\mtf}^{-1}(F_i)\in \cecc$ to be the {\it conjugacy class sums} corresponding to the  conjugacy class $\mathcal C^i$ where  $\lam\in \cfcc$ is  a cointegral such that $\langle \lam, u\rangle =1$.

We denote by $\mu_j:\cfcc\ra \comp$ the algebra homomorphisms corresponding to the central primitive idempotent $F_j$.

We adapt the following from \cite[Definition 1.4]{b-pa}.
\bn{defn}\label{v-property}
A fusion category $\cc$ has \emph{Burnside's vanishing property} (or shortly, is Burnside) if for all $X_i\in  \irrcc $, the following tow assertions are equivalent:  
\bne
\item There is some $\muj\in \wdb$ such that 
$\muj(\ch_{ _{X_i}})=0$,
\item
$X_i$ is not an invertible object.
\ene
\end{defn}

By \cite[Theorem 2]{b-galois} any weakly integral fusion category is Burnside.
\subsection{Definition of the support $\mtc J_{\cd}$}
Let $\cc$ be a fusion subcategory with a commutative Grothendieck ring. For any fusion subcategory $\cd\subseteq \cc$, recall the definition of its support form \cite[Subsection 4.2]{ccc-march}. It is well known that there is an inclusion of $\comp$-algebras $\iota:\cfcd\hookrightarrow \cfcc$.  Thus there is a subset $\mtc J_\cd \subseteq \{0, \dots, m\}$ such that
\beq\label{idemptint}
\lam_{\cd}=\sum_{j \in \mtc J_\cd}F_j
\eeq
 since $\lam_\cd$ is an idempotent element inside $\cfcc$. The set $\mtc J_{\cd}$ is called the \emph{support } of the fusion subcategory $\cd$.
\subsection{The dual hypergroup of $\kc$}
Let $\cc$ be a fusion category with a commutative Grothendieck group $\kc$. On the dual space $\wkc:=\hm_\comp(\kc, \comp)$ one defines a multiplicative structure as in \cite{b-blms}, such as, for all $f, g\in K(\cc)$ one has:
\beq\label{mwa}
[f\star g](\nchi):=f(\nchi)g(\nchi), \;\text{for all}\;\chi_i\in \irrcc.
\eeq

Then $\wkc$ is not anymore a fusion ring but an abelian normalized hypergroup, see \cite[Definition 1.1]{b-pa}.

%
Let $\hb$ be an abelian  normalized  hypergroup. By the commutativity  assumption on $\kc$, the set of algebra morphisms $\widehat{\irrcc}:=\{\muj\;\}_{j=0}^m$ forms a basis for the dual $\wkc$.,

\subsection{The isomorphism $\al$}
Note that the dual hypergroup $\wkc$ is denoted by $\wcfcc$ in \cite{b-blms}. Theorem 3.4 from the same paper implies that for any pseudo-unitary fusion category there is a  canonical isomorphism of hypergroups 
\beq\label{al:isom}
\al:\wcfcc \ra \cecc,\;\muj\sent \frac{C_j}{\dim(\cc^j)}.
\eeq

\section{Proof of Theorem \ref{harada}}\label{pf-harada}
\bpf
Since $\cc$ is weakly-integral, it is a Burnside category, by \cite[Theorem 2]{b-galois}.  Then \cite[Corollary 4.5]{b-pa} implies that 
\beq\label{first}
\big(\prod_{j=0}^m\mu_j\big)^2=\sum_{[X_i]\in \irr(\ccpt)}\wte_i.
\eeq
where $\wte_i$ is the primitive central idempotent of $\wkc$ associated to $\ch_i:=\ch_{X_i}$. Note that $\al(\wte_i)=E_i$, the primitive idempotent of $\cecc$ associated to the irreducible character $E_i$, see \cite[Subsection 6.2]{scalg}. Applying the isomorphism $\al$ to Equation \eqref{first} 
one obtains the following identity in $\cecc$.
\beqn
(\prod_{j=0}^m\frac{C_j}{\fp(\cc^j)})^2=\sum_{[X_i]\in \irr(\ccpt)}E_i.
\eeqn
On the other hand one has by \cite[Eq. (4.13) and Eq. (4.15)]{ccc-march} one has
$$
\frac{\fpcc}{\fp(\ccpt)}\big(\sum_{[X_i] \in \irr(\ccpt)} E_i\big)=\sum_{j \in \mtc J_\ccpt} C_j
$$
and this implies that 
$$
(\prod_{j=0}^m\frac{C_j}{\fp(\cc^j)})^2=\frac{\dim(\ccpt)}{\dimcc}\big(\sum_{j \in \mtc J_\ccpt} C_j\big).
$$
The proof of the theorem  is finished.
\epf
\bc
With the above notations, if the structure constants are integers then one has that 
\beq
\frac{\dimcc}{\dim(\ccpt)}\big| \prod_{j=0}^m \fp(\cc^j)
\eeq
\ec
\bc
If $\cc$ is Burnside and modular with integral structure constants then 
\beq
\fpccad=\frac{\fpcc}{\fp(\ccpt)}\big| \prod_{j=0}^m \fp(d_j)^2
\eeq
\ec
\bpf
It follows that $\fp(\cc^j)=d_j^2$ since $\cc$ is modular. 
\epf

For a semisimple Hopf algebra $H$ recall the notion of commutator left coideal subalgebra $H'$, see \cite{iop}.  In this case Harada's identity can be written as:
\bc
For a semisimple Hopf algebra one has 
$$
(\prod_{j=0}^m \frac{C_j}{\fp(\cc^j)})^2=\blam_{H'}=\frac{1}{\dim(H')}\big(\sum_{\{j\;|\;\cc^j\subseteq H'\}}C_j\big)
$$
\ec
\section{Proof of Theorem \ref{div} and some other applications}\label{pf-mtc}
Suppose that $\cc$ is a weakly-integral modular category. As usually we suppose that $\cc$ is of type $(d_1,n_1;d_2,n_2;\dots d_r,n_r)$ with $d_1$=1. We will use the notation $\fpcc=\fpccone \fpcctw$ where $\fpccone$ is the free square part of $\fpcc$. Thus $\fpccone$ is the unique square free divisor of $\fpcc$ such that $g.c.d(\fpccone, \fpcctw)=1$. Note that $\fpccone\big| n_1$ by \cite[Corollary  1.17]{b-pa}.

Let $\irrcc$ denote the set of isomorphism classes of simple object of $\cc$. For any subset $\ca\subseteq \irr(\cc)$ we denote $\fp(\ca):=\sum_{X\in \ca}\fp(X)^2$.

\br\label{coprime}
if $\cc$ is an integral modular category then  for any $d_i>1$ one  clearly has $\gcd(d_i, \fpccone)=1$ since $d_i^2\big|\fpcc$ by \cite[Theorem 2.11]{eno-weakly}.
\er
\subsection{Proof of Theorem \ref{div}}
Note that Theorem \ref{div} can be reformulated as
\bt 
In any modular category $\cc$ one has $\fpccone\big| n_i$ for all $i\geq 1$.
\et
\bpf
Recall that for any $a>0$ by $\Irr_a(\cc)$ is denoted the set of isomorphism classes of simple objects of $\cc$ of dimension $a$.

Let $G(\cc)$ be the set of invertible objects of $\cc$. Consider the action of $G(\cc)$ on $\irrcc$ and  note that each set $\irr_{d_i}(\cc)$ is a union of orbits of this action. By \cite[Remark 8.17]{eno-annals} one has that $|G(\cc)|\big| \fp(\irr_{d_i}(\cc))=n_id_i^2$.

Thus $\fpccone\big| n_id_i^2$ and since $(d_i, \fpccone)=1$ for all $i$ with $d_i>1$ it follows that $\fpccone \big| n_i$.
\epf

\bp
The first divisibility $d \big| n_1$ is true if $(d, \fp(V)^2)=1$ for any simple $V$ as in \cite{o-yu}.
\ep
\bpf
By \cite[Theorem 1.15]{b-pa} one has that 
$\mtc V(\cc)=\mtc V(\ccpt) \cup \mtc V(d_i^2)$. On the other hand note that
if $p\big| d$ since $(p,d_i^2)=1$ it follows that $p\big| \fpccpt$.
\epf

\br
Note that by \cite[Thm 7. 4]{nat-2014} any modular fusion category $\cc$ with $\fpcc=p^aq^bd$ with $d$ a square free number is WGT. Thus by \cite{o-yu}[Theorem 4.2] these fusion categories decompose as a Deligne product as above.
\er
\subsection{Proof of Theorem \ref{pqr}}
Our proof goes as follows. By \cite[Theorem 4.2]{o-yu}, it suffices to show that $\cc$ is weakly group-theoretical. If $d\neq 1$ then by Theorem \ref{div} it follows that $d\mid n_i$ for all $i>1$. Thus $n_i=dn'_i$ and Frobenius-Perron dimension of $\cc$ can be written as
$$(pqr)^2=\sumitor n'_id_i^2.$$

Then one can do something similar as in the proof of from \cite[Lemma 9.3]{eno-weakly} to deduce that $\cc$ has a symmetric fusion subcategory. If the subcategory $\cc_{pt}$ is degenerate, then its Mueger center is a symmetric fusion subcategory. Otherwise $\cc\simeq \cc_{pt} \boxtimes \cc_{ad}$ and one can repeat the argument for $\cc_{ad}$.
\subsection{Nilpotent fusion categories}
Note that by Theorem 1.11 and Theorem 1.12 of \cite{b-pa} all nilpotent fusion categories with $K(\cc)$ commutative are both Burnside and dual Burnside. Moreover  they are weakly-integral as iterated extensions of $\vect$.
We also prove the following result for nilpotent fusion categories:
\bp\label{dim:nilp}
If $\cc$ is a non-pointed  integral nilpotent fusion category with $K(\cc)$ commutative  then 
\bne
\item 
$\fpccad$ has no square free parts.
\item 
Moreover 
$\fpccone\big|U(\cc)|$ where $U(\cc)$ is the universal grading group of the fusion category $\cc$.
\ene
\ep
\bpf
By \cite[Theorem 1.15]{b-pa} one has $\mtc V(\ccad)=\cup_{i=0}^m \mtc V(d_i)$. For any prime $p$ it follows that if $p\big| \fpccad$ then $p\big| d_i$. By \cite[Corollary 5.3]{NG} then $p^2\big| d_i^2\big| \fpccad$ which shows that $\fpccad$ has no square free part. On the other hand, since $\fpcc=|U(\cc)|\fpccad$ it follows that $\fpccone\big| |U(\cc)|$.
\epf
\bibliographystyle{alpha}
\bibliography{24nov}
\ed
\blue{Rem 8.17 is true for any weakly-integral fusion category! This implies that if the first divisibility is true for any weakly integral fusion category then also the other divisibilities are true for any other weakly-integral fusion categories!}
\br
By \cite[Rem 8.17]{eno-annals} one has $\fpccpt\big| n_id_i^2$ for all $i$. On the other hand if 
Thus if the first divisibility $d\big| n_1$  holds then also the other divisibilities are true for any other weakly-integral fusion categories
 implies that if the first divisibility is true for any weakly integral fusion category !
\er
\blue{Do also something about divisibility in slightly-degenerate fusion categories!}
\blue{For the nilpotent slightly degenerated from YU's paper apply also my result!}
We call an integer \blue{modular-nilpotent}  if every integral modular fusion category of this dimension is nilpotent. 

\newpage
\section{Things to do}
\section{Braided and nilpotent fusion categories}
For a braided category this is nilpotent if and only if the center is also nilpotent.

Suppose that $\cc$ is strictly-weakly integral then $\fpcc^2=$ and it has a group-like element by the freeness.
\subsection{Nilpotent of $p^2d$, $p^3d$ and $p^nd$.}
\bc
Nilpotent of dimension $p^2q$ are pointed. 
\ec
\bc
Let $\cc$ be an integral nilpotent fusion category of square free dimension. Then $\cc$ is pointed. 
\ec
\bpf
Clear from \cite[Corollary 5.3]{NG}.
\epf
\bc\label{nilp:p2d}
Suppose that $\cc$ is a nilpotent fusion category \green{with $K(\cc)$ commutative}  such that $\fpcc=p^2d$. Then $\cc$ is pointed.
\ec
\bpf
By the previous Theorem one has $d\big| |U(\cc)|$. Since $\ccad$ has no square-free part it follows that $\fpccad=p^2$. Then from FPdim equation one has $p^2\big| |G(\cc)|$. Thus 

$d=k_11^2+n_2$ which implies that $k_1\big| n_2$.

$|G(\cc)|=p^2k_1$ it follows that there are $k_1$ components containing invertible objects. It follows that all the other $d-k_1=n_2$ contains a single simple object of dimension $p$.

$|U(\cc)|=d$ is abelian and therefore products of cyclic groups of prime order.

The set of components containing only invertible objects determine a subgroup $U_1$ of $U(\cc)$.

{\bf case 1: only one prime $d=q$ and $ \fp(\cc)=p^2q$!}

In this case $U(\cc)$ is simple so the fusion category is pointed.

{\bf case 2: $d=q_1q_2$, $\fp(\cc)=p^2q_1q_2$}

Suppose that $U_1=\mathbb Z_{q_1}$ and $U_\cc=<a>\times <b>$. Then

$|G(\cc)|=q_1p^2$ is abelian and therefore $\gcc\simeq G(\ccad)\times \mathbb Z_{q_1}$.

The objects from $\ccad$ fix any $X_b$ with $\fp(X)=p$. Thus in this case the other invertible objects do not fix $X$. Thus the orbits have dimension $q_1$. 

For $g\in \mathbb Z_{q_1}$ either fixes a point or the orbit has dimension $p$. 

If there is no fixed point then there are $q_2-1$ objects of dimension $p$. Thus one needs to have $q_1\big| q_2$.

Then $G(\ccad)$ and $\{X_b\;|\fp(X_b)=p\}$ form a fusion subcategory of dimension 

$p^2+$

and the dimension Equation can be written as
$$p^2q_1q_2=q_1p^2+q_1(q_2-1)p^2$$

Let $X_b\in \cc_b$.  Since $X_b^2\in \cc_{b^2}$
which shows that $$X_b^2=pX_{b^2}$$

One has $X_{b^*}\cc_{b^{-1}}$ so in the odd case there are no self-dual objects.

Suppose 
$$\cc_b\subseteq \ccpt$$

$$XX^*=\sum_{g \in \ccad}g$$

Thus any componente of the universal grading has $\fp$ equal to $p^2$  which shows that it contains either one simple of dimension $p$ or only grouplike elements. 

This is impossible for $\ccad$. 
\epf
\bc
If ??
\ec
\bc
Nilpotent fusion categories of dimension $2p^2d$ are pointed.
\ec
\bpf
If it has a self-dual object of dimension $p$ it should be from $\mathbb Z_2$.

If it has no self-dual objects then $X_b^r=p^{r-1}X_{b^r}$. 
\epf
\bc\label{nilp:p3d}
Suppose that $\cc$ is a nilpotent fusion category such that $\fpcc=p^3d$. Then $\cc$ is pointed.
\ec
\green{Do not identify $\ccpt$ with $U_\cc$ .
}
\bpf
By the previous Theorem one has $d\big| |U(\cc)|$. Since $\ccad$ has no square-free part it follows that $\fpccad\in\{p^2,p^3\}$.
{\bf case 1} $\fpccad=p^2$. Every non invertible object has dimension $p$. The same argument for $\ccad$ implies that $p^2\big| \fpccpt$ which gives that  
Thus any componente of the universal grading has $\fp$ equal to $p^2$  which shows that it contains either one simple of dimension $p$ or only grouplike elements. This is impossible for $\ccad$. 
\epf
\bibliographystyle{alpha}
\bibliography{24nov}
\newpage
\newcommand{\rctwo}{{R_{\cc^{(2)}}}}
\newcommand{\rg}{{\mtr{rank}(\cc^g)}}
\newcommand{\kcd}{{K(\cd)}}
\newcommand{\wkcd}{\widehat{\kcd}}
\subsection{On gradings and dual Burnside}

\bl
For any graded extension $\cc=\bigoplus_{g \in G}\cc_g$ one has
$$
Pr_1=r_{ _{\prod_{g \in G}}g^{\rank(\cc_g)}}.
$$
where $r_g=\frac{1}{\fp(\cc_1)}R_g$ is the normalized regular character of $\cc_g$.
\el
\green{On the other hand all the constituents of $P$ are inside $\ccad\subseteq \cc_1$. This implies that $\prod_{g \in G}g^{\rank(\cc_g)}=1$.
}
\bpf 

Note that for any $x\in \cc_g$ one has 
$xR_1=\fp(x)R_g$ if $g\in G$.  This implies that $R_gR_h=\fp(R_g)R_{gh}$ for all $g,h\in G$.

One has
$$
P:=\prod_{x\in \cc} \frac{x}{d_x}=\prod_{g \in G}\prod_{x \in \cc_g}\frac{x}{d_x}.
$$
Denote by $P_g:=\prod_{x\in \cc_g}\frac{x}{d_x}$. Therefore  $P=\prod_{g \in G}P_g$. Then 
$$
P_g r_1=P_gr_1^\rg=\prod_{x \in \cc_g}(\frac{x}{d_x}r_1)=\prod_{x \in \cc_g}r_g=r_g^\rg.
$$

From here we obtain 
$$
Pr_1=\prod_{g \in G}(P_gr_1)=\prod_{g \in G}r_{g^\rg}=r_{\prod_{g \in G}g^{\rg}}
$$
\epf
\green{
If $\cc_1$ is dual-Burnside then $P_1=R_{ _{{\cc_1}_{\ad}}}=\rctwo$.}

\bne
\item 
Note that $P_gR_g=P_gR_1^{\rg}=R_g^{\rg}$
\item $P=R_{ _{\prod_{g \in G} g^{\rg}}}$
\item
Then  $P_g\rctwo$
\item suppose that $\cc^2=\cc^1$ (i.e the dual is perfect).
\item Then $\rctwo=R_{\cc_1}$. Therefore one has
$$
P=\prod_{g \in G}P_g=
$$
\ene
\blue{It is related to the Hecke algebra of $\cc^{2}$ inside $\cc$.
\\
Suppose that $\cc_1$ is pointed. Then $r_gr_h=r_{gh}$
\\
Suppose that $\cc^{(2)}\subseteq \cc^{(1)}\subseteq \cc$
}
\bt
If $\cc=\oplus_{g\in G}\cc^g$ and $\cc_1$ is dual-Burnside then $G$ is the universal grading group and $$\prod_{g \in G}g^{\rg}=1.$$ 
and $\cc$ is also dual-Burnside.
\et
\bpf
By the previous lemma one has

$$P=R_{\prod_{g \in G}}g^{\rank(\cc_g)}$$
\blue{Suppose now that $\cc$ is a dual-Burnside. Then $P$ has to be equal to $R_{\ccad}$. On the other hand }

Suppose now that $\cd:=\cc_1$ is dual Burnside.

On the other hand the constituents of $P$ are all inside $\ccad\subseteq \cc_1$. It follows that $\prod_{g \in G}g^{\rg}=1$ and one has the universal grading, i.e $\cc_1=\ccad$.
\epf
\bl
If $\cc$ is not nilpotent then there is a situation with 
$$
\cc\supseteq {\cc_{ad}}={(\cc_{ad})}_{ad}.
$$
\el
\ed